\documentclass[12pt]{amsart}
\numberwithin{equation}{section}
\newtheorem{theorem}{Theorem}

\newtheorem{proposition}{Proposition}
\numberwithin{theorem}{section}
\numberwithin{lemma}{section}
\numberwithin{proposition}{section}
\newtheorem{corollary}{Corollary}
\newcommand{\ignore}[1]{}
\def \reel{ {\rm I}\!{\rm R} }
\begin{document}
\tracingpages
1
\title[heat kernel]{Large time behavior of heat kernels on forms}
\author{Thierry Coulhon and Qi S. Zhang}
\address{ D\'epartement
de Math\'ematiques, Universit\'e de
Cergy-Pontoise, 95302 Pontoise, France
}
\address{ Department of Mathematics, University of
California, Riverside,
CA 92521, USA }
\date{\today}
\thanks{TC's research was partially supported by the European
Commission
(IHP Network ``Harmonic Analysis and Related Problems''
2002-2006, Contract HPRN-CT-2001-00273-HARP)}

\begin{abstract}
We derive large time upper bounds for heat kernels on
vector bundles of differential forms on a class of non-compact
Riemannian manifolds under  certain curvature conditions. 
\end{abstract}
\maketitle
\tableofcontents
\section{Introduction}

The
goal of the present paper is to establish large time, pointwise
bounds for the heat
kernel on the vector bundles of forms on some
noncompact manifolds.  

Information on large time behavior of
heat kernels on forms
usually leads to interesting analytical and
topological
information on the manifolds. In fact heat kernel on
forms
contains much more information on the interplay between
analysis,
geometry and topology than that on functions. So far much
effort
has been spent on the study of short time and long time behavior of heat kernel
on
forms, in the  case of closed manifolds, see for instance \cite{BGV}, \cite{Ro}. By
contrast, the present paper is to our knowledge the first one to offer estimates 
for the heat kernel on one-forms on a class of non-compact Riemannian manifolds
with a meaningful contents for large time, i.e. without an increasing exponential factor
(see for instance \cite{Lo}).

Let $M$ be a complete connected Riemannian manifold. Denote by $d(x,y)$ the geodesic
distance between two points
$x,y\in M$,  and by $B(x,r)$ the open ball of center $x\in M$. Let
$\mu$ be the Riemannian measure; denote also by
$|\Omega|$ the measure $\mu(\Omega)$ of a mesurable subset
$\Omega$ of $M$.  Denote by $\Delta$ the (non-negative)
Laplace-Beltrami operator on functions. The heat semigroup
on functions $e^{-t\Delta}$  will also be denoted by
$P_t$, and the corresponding heat kernel by $p_t(x,y)$, $t>0$, $x,y\in M$.
  We will
use $\vec{\Delta}$ to denote the Hodge
Laplacian on
  forms. The heat semigroup
on forms $e^{-t\vec{\Delta}}$  will also be denoted by
$\vec{P}_t$, and the corresponding heat kernel by $\vec{p}_t(x,y)$, $t>0$, $x,y\in M$.

The main
question we shall address below is the following:

{\it Given an upper estimate for  the heat kernel on
functions, under which additional assumptions can one deduce an upper
bound for the heat
kernel on forms?}

We shall consider in particular the case where $M$ has the so-called volume doubling property and the
heat kernel on functions satisfies a Gaussian upper estimate, that is
\[
p_t(x, y) \le \frac{C}{|B(x,
\sqrt{t})|} \exp(-c d^2(x, y)/t), \
\forall\,x,y\in M,\, t>0,
\]
for some
$C,c>0$.

For instance, when $M$ has non-negative 
Ricci curvature, it was proved in \cite{LY}
that $p_t$ satisfies a Gaussian upper estimate, and in that case,   the answer to the above
question for $1$-forms is
straightforward by the semigroup domination theory
(see  \eqref{dom} below, and also e.g. \cite{HSU}, \cite{HSU2}, 
\cite{DL}, \cite{R}): the heat kernel on $1$-forms is also
bounded from above by a
Gaussian.

 The following simple example shows that this may be false in general.  Let $M$ be the connected sum of two
copies
of $\reel^n$, $n \ge 3$. It is known that the heat kernel
on functions
has a Gaussian upper bound (see \cite{BCF}).
If the heat kernel on 1-forms $\vec{p}_t$ also had a Gaussian
upper
bound, then by \cite{CD}, pp.1740--1741, it would follow
that \[
|\nabla p_t(x, y)| \le \frac{C'}{\sqrt{t} |B(x,
\sqrt{t})|} \exp(-c' d^2(x,
y)/t),\ \forall\,x,y\in M,\, t>0.
\]
A classical argument shows that $p_t$ would then be  bounded
below by a Gaussian (see for instance \cite{LY}). This is false as was noticed in
\cite{BCF}. See   \cite{CD0} for more on this example.

Another case where the behaviour of the heat kernel on $1$-forms  is well understood is the case of the Heisenberg
group and more generally stratified Lie groups, see \cite{R1}, \cite{R2}.

In the present  paper, we are going to see that if the negative part of the Ricci curvature is small enough in some sense,
then  the upper
bound on the heat kernel on $1$-forms differs from that on the heat
kernel
on functions at most by a certain power of time $t$. One can state similar results  for higher degree forms by replacing in the assumptions
the Ricci curvature by a suitable curvature operator (see \cite{GM}).
For convenience, we shall however formulate our assumptions and results in the case of $1$-forms.
We leave the formulation of the general case to the reader.

\bigskip

In this article, all Riemannian manifolds under consideration will be complete non-compact.
Let us layout some basic assumptions to
be used below.

\medskip

{\it Assumption (A). $M$ satisfies the volume doubling
property.
$$
|B(x,
2r)| \le C |B(x,
r)|
$$
for all $x \in M$, $r>0$
and some $C> 0$.}

\medskip

{\it Assumption (B). The heat kernel  $p_t(x, y)$ on functions
satisfies a
Gaussian upper bound:
$$
p_t(x, y) \le \frac{C}{|B(x,
\sqrt{t})|} \exp(-c d^2(x,
y)/t),
$$
for  some $C, c>0$, and  all $x, y \in  M$ and $t>0$.}

\medskip

It
was proved in \cite{G:1} that
Assumptions $(A)$ and $(B)$ together are
equivalent to the
following relative Faber-Krahn inequality:

For all $x\in M$, $r>0$, and every non-empty subset $\Omega \subset B(x,
r)$,
\[
\lambda_1(\Omega) \ge \frac{c}{r^2} \left( \frac{|B(x,
r)|
}{|\Omega|}
\right)^{2/\nu}.
\leqno{(FK)}\]Here $\lambda_1(\Omega)$ is the first Dirichlet
eigenvalue of $\Omega$ and $c>0$.

Note that $(FK)$ implies
\begin{equation}
\frac{|B(x,s)|}{|B(x,r)|}\le C\left(\frac{s}{r}\right)^\nu,\label{dou}
\end{equation}
for all $s>r>0$, $x\in M$, and we shall use Assumption $(A)$ in this form.

\medskip

{\it
Assumption (C). The Ricci
curvature is
bounded from below by a negative constant.}

\medskip

It follows from Assumption $(C)$ and Bishop's comparison theorem that there exists $C>0$ such that $ |B(x,
1)| \le C$ for all $x \in  M$. We shall often also need the opposite inequality.

\medskip

{\it
Assumption (D). Non-collapsing  of the volume of balls: there exists $c>0$ such that $ |B(x,
1)| \ge c$ for all $x \in  M$.}

\bigskip

It is well-known that to estimate the heat kernel acting on one-forms, it is enough to estimate the  kernel
 of a certain  Schr\"odinger semigroup acting on functions, whose potential is the negative part of the Ricci
curvature.

  Indeed, Bochner's formula states
  $$
\overrightarrow{\Delta} = D^* D - Ric.
  $$
Here $D$ is the
covariant derivative on $1$-forms and $Ric$ is the Ricci curvature.

  Let
$\lambda=\lambda(x)$ be the lowest eigenvalue of $Ric(x)$, $x\in
 M$.
We
  will use the notation
$$
  V(x) = \lambda^-(x) = (
|\lambda(x)|-\lambda(x) ) /2.
$$

Let $P_t^V$ be the semigroup $e^{-t(\Delta - V)}$.  Under Assumption $(C)$,   $P_t^V$ has a  kernel which we shall denote by $p_t^V(x, y)$.
 
 Let us recall the semi-group domination
property, which  was proved in \cite{HSU2}:
\begin{equation}
|\vec{p}_t(x, y) | \le  p_t^V(x, y)
\label{dom}
\end{equation}
for all $x, y \in  M$ and $t>0$.  Here
$\vec{p}_t(x, y)$ is a linear operator from  the tangent space $T_yM$ to $T_xM$, and
here $|\vec{p}_t(x, y) |$ denotes its operator norm with respect to the Riemannian  metrics.

\bigskip

We can now introduce one of our main curvature assumptions.
An
important property of  the Hodge Laplacian
$\overrightarrow{\Delta}$ is that
it is a nonnegative operator (as a consequence, $\vec{P}_t$ is contractive on $L^2(M,T^*M)$). This means that, for every  smooth
compactly supported
$1$-form $\phi$, 
\begin{equation}
\int_ M -Ric( \phi(x),
\phi(x)) \,d\mu(x) \le  \int_ M | D
\phi (x)|^2 \,d\mu(x).
\label{ri}
\end{equation}
Now, by the Kato inequality
$$
|\nabla |\phi||\le | D
\phi |,
$$
and the fact that by definition
$$Ric( \phi(x),
\phi(x)) \ge -V(x)|\phi (x)|^2,
$$
we see that condition (\ref{ri})
is implied by
\[
\int_{M} V(x) f^2(x) \,d\mu(x) \le  \int_{M} | \nabla
f(x)
|^2 \,d\mu(x),\  \forall\,f \in C^{\infty}_0( M),
\]
which means that $\Delta -V$ is a positive operator on $L^2(M)$.

We shall say that  $\Delta -V$ is strongly positive (strongly subcritical in the sense of \cite{DS}) if it satisfies the following stronger condition:  there exists
$A<1$ such that, for all $f \in C^{\infty}_0( M)$,
\begin{equation}
\int_{M} V f^2 \,d\mu \le A \int_{M} | \nabla
f
|^2 \,d\mu.
\label{sp}
\end{equation}

 The
above condition sometimes is
referred to as the form boundedness condition,
which has its
origin in the Hardy type inequality, for $f \in
C^\infty_0(\reel^n)$, $n \ge 3$,
\[
\frac{(n-2)^2}{4} \int_{\reel^n}
\frac{f^2(x)}{|x|^2}  dx \le
\int_{\reel^n} |\nabla f(x) |^2 dx.
\]
For generalizations of the above inequality to the
manifold case, see \cite{Ca}. 

\medskip

{\it  Example:}  If  the manifold $M$ satisfies the
Euclidean Sobolev inequality of dimension $n$
$$
\left(\int_{M}  |f|^{2n/(n-2)} \,d\mu\right)^\frac{n-2}{n} \le C \int_{M} | \nabla
f
|^2 \,d\mu
$$
for all $f \in
C^\infty_0(M)$, for some $n>2$, and if 
$V\in L^{n/2}(M)$ with sufficiently small norm, then it is easy to see by using H\"older's inequality that (\ref{sp}) holds.

\bigskip

Let us now summarize our  results. Under Assumptions $(A)$ to $(D)$, the function $V$, the negative part of the lowest eigenvalue
of the Ricci curvature largely determines the upper bound of heat hernel 
on $1$-forms. If  $V$ is sufficiently small in certain integral sense, 
then $\vec{p}_t$ has Gaussian upper bound, which has important consequences in terms of $L^p$ boundedness of the Riesz transform. This is explained in Section \ref{rie}.
Otherwise the upper bound for $\vec{p_t}$ is a Gaussian times a suitable power of time
$t$, provided that the operator $\Delta - V$ is strongly positive. The proof of this fact is contained in Sections
\ref{poin}, \ref{1}, \ref{end}.  We also consider the case where the Ricci curvature is nonnegative outside of a compact set. Without any other assumptions on the Ricci curvature, there may be $L^2$ harmonic forms, therefore one cannot expect a  decay with respect to time in general, but we show that $\vec{p_t}$ is bounded by a Gaussian plus the product of the Green's function of the Laplacian in both variables. This is the subject of Section \ref{bad}.
Finally, we treat  in Section \ref{wd} the case where the heat kernel on functions has an arbitrary uniform decay.

\maketitle
\section{Bounds for the heat kernel on forms and strong positivity}\label{stro}

Our aim in this section is the following result.

\begin{theorem}
\label{main}
Suppose $M$
satisfies Assumptions $(A)$, $(B)$, $(C)$, and $(D)$, and  that the
operator $\Delta -V$ is strongly positive with constant $A$.
Suppose in addition that $V \in L^p( M,\mu)$ for some $p\in[1,+\infty)$. Then, 
 if $p=1$,  for any $0<c<1/4$,  and any $\varepsilon>0$, there exists
$C>0$ such that
 \[
|\vec{p}_t(x, y)| \le C \min \left\{ \frac{ 
t^{ (1+\varepsilon)A}}{|B(x,
\sqrt{t})|}, \, 1 \right\} \exp (-
c d^2(x, y)/t),  \ \forall\, x, y \in  M, \ t \ge 1;\]
if $p \in (1, 2)$,  for any $c<1/4$, there exists
$C, \varepsilon>0$ such that
 \[
|\vec{p}_t(x, y)| \le C \min \left\{ \frac{ 
t^{(p-\varepsilon) A}}{|B(x,
\sqrt{t})|}, \, 1 \right\} \exp (-
c d^2(x, y)/t),  \ \forall\, x, y \in  M, \ t \ge 1;
\]If $p \ge 2$, 
for any $c<1/4$ and any $\varepsilon>0$, there exists $C$ such
that
 \[
|\vec{p}_t(x, y)| \le C \min \left\{ \frac{ 
t^{(p-1+\varepsilon) A}}{|B(x,\sqrt{t})|}, \, 1 \right\}
\exp (- c d^2(x, y)/t),  \ \forall\, x, y \in  M, \ t \ge 1.
\]
\end{theorem}

{\it Remark:} 
Using Assumption $(C)$ and the Gaussian bound on $p_t$, one easily obtains the following small time estimate 
\begin{equation}
p_t^V(x, y) \le  \frac{C_{t_0}}{|B(x,
\sqrt{t})|}\exp (- c d^2(x, y)/t),
\label{small}
\end{equation}
 for all $x,y\in M$ and $0<t<t_0$.

Thanks to the domination property \eqref{dom}, Theorem \ref{main}
is a consequence of the following statement, which is of independent interest.

\begin{theorem}
\label{mainV}
Suppose $M$
satisfies Assumption $(C)$, and  that the
operator $\Delta -V$ is strongly positive with constant $A$.
Suppose in addition that $V \in L^p( M,\mu)$  for some $p\in[1,+\infty)$. Then, 
 if $p=1$,  for any $0<c<1/4$ and any $\varepsilon>0$, there exists
$C>0$ such that
 \[
p^V_t(x, y) \le C \min \left\{ \frac{ 
t^{ (1+\varepsilon)A}}{|B(x,
\sqrt{t})|}, \, 1 \right\} \exp (-
c d^2(x, y)/t),  \ \forall\, x, y \in  M, \ t \ge 1;\]
if $p \in (1, 2)$,  for any $c<1/4$, there exists
$C, \varepsilon>0$ such that
 \[
p^V_t(x, y) \le C \min \left\{ \frac{ 
t^{(p-\varepsilon) A}}{|B(x,
\sqrt{t})|}, \, 1 \right\} \exp (-
c d^2(x, y)/t),  \ \forall\, x, y \in  M, \ t \ge 1;
\]If $p \ge 2$, 
for any $c<1/4$ and any $\varepsilon>0$, there exists $C$ such
that
 \[
p^V_t(x, y) \le C \min \left\{ \frac{ 
t^{(p-1+\varepsilon) A}}{|B(x,\sqrt{t})|}, \, 1 \right\}
\exp (- c d^2(x, y)/t),  \ \forall\, x, y \in  M, \ t \ge 1.
\]
\end{theorem}

We would like to mention a number of previous papers that deal with Schr\"odinger heat kernels on manifolds. In the paper \cite{LY},  a fundamental gradient estimate was derived
for the heat kernel. As far as long time behavior is concerned, the emphasis is on the case without potential and nonnegative Ricci
curvature.  The paper \cite{S} studied Schr\"odinger heat kernels with singular oscillating potentials.  The papers  \cite{Z1,Z2}  established long time behavior for Schr\"odinger heat kernels
on manifolds with nonnegative Ricci curvature for potentials essentially behaving as negative powers
of the distance function. The case of potentials with polynomial growth and magnetic field is considered in \cite{Si3}.

Here is the plan of the proof of Theorem \ref{mainV}.

 In section \ref{poin},  we show that, given the upper bound on $p_t$
and the strong positivity of $\Delta-V$, a pointwise upper bound on $p_t^V$
follows from an adaptation of the Nash method due to Grigor'yan, provided some $L^1$ to $L^1$ estimates for $P_t^V$ are available.
In section \ref{1}, we prove such estimates under the other assumptions of Theorem \ref{mainV}, and we finish the proof.

\maketitle
\subsection{Pointwise estimates}\label{poin}

Let us first prove the   following preliminary estimate.

\begin{proposition}\label{prel}
Suppose $M$
satisfies Assumption $(C)$, and  that the
operator $\Delta -V$ is strongly positive. 
Then there exist $C,c>0$ such that
\begin{equation}
p_t^V(x, y)  \le C \exp\left(-c d^2(x, y)/t\right), \ \forall \,x,y\in M,
\ t\ge
1.\label{exp}
\end{equation}

\end{proposition}


\proof

This can be proven by a standard
method of using exponential
weights as in \cite{Gmax}. An alternative way is to use wave
equation method as in
\cite{CGT} or \cite{Sik}.  

Fix $y$ and write
\[
u(x, t) =   p^V_t(x, y),\]
\[
I(t) = \int_M  u^2(x, t) w(x, t) d\mu(x)
\]where $w(x, t) = e^{ \frac{d^2(x, y)}{D t}}$ for some $D>0$ to be chosen
later.

One has
\begin{eqnarray*}
 \frac{d}{dt} I(t) &= &\frac{d}{dt} \int_M  u^2 w \,d\mu\\
 &=&2\int_M   u w (- \Delta  u + V u) \,,d\mu
 - \int_M   \frac{u^2 w d^2(x, y)}{D t^2} d\mu(x).
\end{eqnarray*}
This implies, after
integration by parts,
\begin{eqnarray*}
\frac{d}{dt} I(t)
& = &
 - 2\int_M   |\nabla u|^2 w  \,d\mu  -2\int_M  u\nabla w\cdot\nabla u \,d\mu\\
 &&+
 2\int_M   V (u \sqrt{w})^2  \,d\mu
 - \int_M   \frac{u^2 w d^2(x, y)}{D t^2} d\mu(x)\\
 & = &
 - 2\int_M   |\nabla u|^2 w  \,d\mu  -2\int_M  uw\frac{2 d(x, y}{D t})\nabla d(x, y)\cdot\nabla u\, d\mu(x)\\ &&
 +2\int_M   V (u \sqrt{w})^2  \,d\mu
 - \int_M   \frac{u^2 w d^2(x, y)}{D t^2} d\mu(x)\\
 & \le &
 - 2\int_M   |\nabla u|^2 w  \,d\mu  +4\int_M  uw\frac{ d(x, y)}{D t}|\nabla u | d\mu(x)\\
&&+ 2\int_M   V (u \sqrt{w})^2  \,d\mu
 - \int_M   \frac{u^2 w d^2(x, y)}{D t^2} d\mu(x)\\
& \le &
 - 2\int_M   |\nabla u|^2 w  \,d\mu +C_\varepsilon\int_M  \frac{ u^2wd^2(x, y)}{D^2 t^2}d\mu(x)\\
 &&+\varepsilon\int_M  |\nabla u |^2 w\,d\mu+
2\int_M   V (u \sqrt{w})^2  \,d\mu
 - \int_M   \frac{u^2 w d^2(x, y)}{D t^2} d\mu(x)\\
 & \le &
 - (2-\varepsilon)\int_M   |\nabla u|^2 w  \,d\mu +
2\int_M   V (u \sqrt{w})^2  \,d\mu
 - \int_M   \frac{u^2 w d^2(x, y)}{2D t^2} d\mu(x)
 \end{eqnarray*}
 for arbitrarily small $\varepsilon>0$,
 provided $D$ is chosen large enough.
 Finally, using the strong positivity of $V$, we obtain
$$\frac{d}{dt} I(t)\le  - (2-\varepsilon)\int_M   |\nabla u|^2 w  \,d\mu  +
 2 A \int_M   |\nabla (u \sqrt{w})|^2  \,d\mu
 - \int_M   \frac{u^2 w d^2(x, y)}{2D t^2} d\mu(x).$$

Observe that
\begin{eqnarray*}
|\nabla (u \sqrt{w})|^2 & =& |\nabla u |^2 w + u \nabla u \cdot\nabla w +
\frac{u^2 |\nabla w |^2}{4 w}\\
& \le& |\nabla u |^2 w +  2 u |\nabla u| \frac{d(x, y)}{D t} w +
\frac{d^2(x, y) u^2 w}{D^2 t^2}\\
& \le& |\nabla u |^2 w +  2C \frac{d^2(x, y)}{D^2 t^2}  u^2 w +  \frac{2}{C} |\nabla u|^2 w +
\frac{d^2(x, y) u^2 w}{D^2 t^2}.
\end{eqnarray*}
  We find, since $A<1$,
that
\[
\frac{d}{dt}I(t) = \frac{d}{dt} \int_M  u^2(x,t) e^{ \frac{d^2(x, y)}{D
t}} d\mu(x) \le
0
\]when $D$ is sufficiently large.

In particular,
\[
\int_M  \left(p^V_t(x, y)\right)^2   e^{ \frac{d^2(x, y)}{D t}} d\mu(x)
\le \int_M   \left(p^V_1(x, y)\right)^2   e^{ \frac{d^2(x, y)}{D
}} d\mu(x)
\]when $t \ge 1$. By
\eqref{small} and the small time Gaussian estimate for $p_t$ under Assumption $(C)$
(see \cite{LY}),
\[
p^V_1(x, y) \le
\frac{C}{|B(x, 1)|} e^{-c d^2(x, y)}.
\]Using the well-known fact that a manifold satisfying $(C)$ has at most exponential volume growth around any point, we have, for $D>0$ large enough,
\[
\int_M  \left(p^V_1(x, y)\right)^2   e^{ \frac{d^2(x, y)}{D }} d\mu(x)
\le C, 
\]

hence
\[
\int_M   \left(p^V_t(x, y)\right)^2   e^{ \frac{d^2(x, y)}{D t}} d\mu(x)
\le C, \quad t \ge 1.
\]

Next, using the semigroup property
\begin{eqnarray*}
 p^V_{2t}(x, y) &=& \int_M  p^V_t(x, z)
p^V_t(z, y) d\mu(z)
\\
&=&\int_M  e^{ \frac{d^2(x, z)}{2 D t}} p^V_t(x, z) e^{
\frac{d^2(z, y)}{2D t}} p^V_t(z, y) e^{-\frac{d^2(x, z)}{2 D
t} -
\frac{d^2(z, y)}{2 D t}}  d\mu(z)\\
&\le&  e^{ -\frac{d^2(x, y)}{4 D t}} \big{[} \int_M  e^{
\frac{d^2(x, z)}{ D t}} \left(p^V_t(x, z)\right)^2 dz \big{]}^{1/2}
\big{[} \int_M  e^{ \frac{d^2(y, z)}{ D t}}\left( p^V_t(y,
z)\right)^2  d\mu(z)\big{]}^{1/2}.
\end{eqnarray*}
Hence
\[
p^V_{t}(x, y) \le C e^{-c d^2(x, y)/t}, \quad t \ge 2.
\]This proves the claim.
\qed

\medskip

We can now state our main technical result, which is an adaptation of the Nash method to the
case where the semigroup under consideration is not necessarily contractive on $L^1$.
The argument is
based on the one  in the proof of  \cite{G}, Theorem 1.1 (see also \cite{CG}, Proposition 8.1), with certain
modification and
localization. If $T$ is a operator from $L^{p_1}$ to $L^{p_2}$,
then $\Vert T
\Vert_{p_1, p_2}$ will denote the operator norm $\sup_{f
\in
L^{p_1}\setminus\{0\}} \frac{ \Vert T f \Vert_{p_2}}{\Vert f \Vert_{p_1}}.$

\begin{proposition}\label{nash}
  Let $M$ satisfy 
Assumptions (A), (B), (C), and (D). Suppose that $\Delta-V$ is strongly
positive
and that
there exists an non-decreasing function $F$ such that
  \begin{equation}\label{l1}
\Vert P_t^V
\Vert_{1,1} \le F(t), \ t\ge 1.
\end{equation}
 Then there exists $C>0$ such that
$$
p^V_t(x, x)\le C \frac{F^2(t) [\ln (e+ t
F(t))]^{\nu/2}}{|B(x,\sqrt{t})|},
$$
for all $x\in M$ and $t\ge 1$.
where $\nu>0$ is the constant from \eqref{dou}.
If in addition $F$ satisfies
$$F(2t)\le CF(t),Ê\forall\,t>0,$$
 for some $C>0$, then
 for any $c\in (0,1/4)$, there exists $C>0$ such that \begin{equation}
p^V_t(x, y)\le C \frac{F^2(t) [\ln (e+ t
F(t))]^{\nu/2}}{|B(x,\sqrt{t})|} \exp (- c d^2(x, y)/t), 
\label{im}
\end{equation}
for all $x,y\in M$ and $t\ge 1$.
\end{proposition}

{\it Remark.} If there exists $N$ such that $F(t)\le Ct^N$, the above estimate takes the following simpler form:
$$
p^V_t(x, y)\le C \frac{F^2(t) [\ln (e+ t
)]^{\nu/2}}{|B(x,\sqrt{t})|} \exp (- c d^2(x, y)/t).
$$

\proof

Fix $x_0\in M$,   write $u(x, t) = p_t^V(x, x_0)$, $t>0$, $x\in  M$, and set
\[
I(t) =
\int_ M u^2(x, t)\,d\mu(x)=p_{2t}^V(x_0, x_0).
\]Then
\[
I'(t) = -2 \int_ M u(x, t)
(\Delta u - V u)(x, t) \,d\mu(x)=-2\int_ M | \nabla u
|^2 \,d\mu+2\int_ M V u^2
\,d\mu 
.
\]
Using  assumption (\ref{sp}), we have 
\begin{equation}
I'(t) \le -2 (1-A) \int_ M |\nabla u(x, t)|^2 \,d\mu(x).
\label{diff}
\end{equation}

Since, for any $s>0$,
\[
u^2 \le (u-s)^2_+ + 2 s u,
\]we
can write
\[
I(t) \le \int_{ \{x | u(x, t) >s \}} (u(x,
t)-s)^2 \,d\mu(x) + 2s
\int_ M u(x, t) \,d\mu(x).
\]
By  assumption
(\ref{l1}) and the definition of $\lambda_1$, this yields
$$
I(t) \le
\frac{\int_{ \{x | u(x, t) >s \}} |\nabla (u(x,
t)-s)|^2
\,d\mu(x)}{\lambda_1 ( \{x | u(x, t) >s \})}  + 2s F(t),
$$
hence
\begin{equation}
I(t) \le
\frac{\int_{ \{x | u(x, t) >s \}} |\nabla u(x,
t)|^2
\,d\mu(x)}{\lambda_1 ( \{x | u(x, t) >s \})}  + 2s F(t).
\label{lev}
\end{equation}

The bound \eqref{exp} yields
\[
u(x, t) \le C_1 e^{- c_2
d^2(x, x_0)/t},
\]thus
\[
  \{x \ | \ u(x, t) >s \} \subset \{ x \ | \  e^{-
c_2 d^2(x, x_0)/t} > s/C_1
\} = \{ x \ | \ d^2(x, x_0)< c^{-1}_2 t  \ln
(C_1/s)  \}.
\]
Thus
\[
  \{x \ | \ u(x, t) >s \} \subset B(x_0,r),
\]
where 
\[
r = \sqrt{c^{-1}_2 t (| \ln
(C_1/s)|+1)}
\]
(we choose to take $r\ge c\sqrt{t}$ for later convenience).
According to $(FK)$, we
have
\[
\lambda_1(\{x \ | \ u(x, t) >s \}) \ge \frac{c}{r^2}
\left(
\frac{|B(x_0,r)|}{|\{x \ | \ u(x, t) >s \}| } \right)^{2/\nu}.
\]
On the other hand,
\[
|\{x \ | \ u(x, t) >s \}| \le s^{-1} \int_ M u(x, t) \,d\mu(x)
\le
s^{-1} F(t).
\]Therefore
\begin{equation}
\lambda_1(\{x \ | \ u(x, t) >s
\}) \ge \frac{c}{r^2} \left(
\frac{s|B(x_0,r)|}{ F(t) } \right)^{2/\nu}
:= m(s, t, x_0).
\label{fab}
\end{equation}

Plugging this into
\eqref{lev}, we obtain
\[
I(t) \le \frac{\int_{ \{x | u(x, t) >s \}} |\nabla
u(x, t)|^2
\,d\mu(x)}{m(s, t, x_0)}  + 2s
F(t).
\]
Hence
\begin{equation}
\int_{ \{x | u(x, t) >s \}} |\nabla u(x,
t)|^2 \,d\mu(x) \ge \left( I(t) - 2 s
F(t) \right) m(s, t, x_0).
\label{sub}
\end{equation}

The combination of  (\ref{sub}) and (\ref{diff})
yields
$$
I'(t) \le - 2(1-A) \left( I(t) - 2 s
F(t) \right) m(s, t, x_0).
$$
Choosing $s$ so that $s F(t) = I(t)/4$ yields
\begin{equation}
I'(t) \le -
(1-A) I(t)  m(s, t, x_0),\label{bec}
\end{equation}
for all $t>0$ and the  corresponding $s$.

We have that
\[
I(t) =p^V_{2t}(x_0, x_0) 
\ge c/t^{\nu/2},
\]
for $t\ge 1$. Indeed, since $V \ge 0$, by the maximum
principle
$p_t^V(x, x_0) \ge p_t(x, x_0)$.
Now it is well known (see \cite{BCF}) that Assumptions $(A)$ and $(B)$ imply
$$p_{2t}(x_0, x_0)\ge \frac{c}{|B(x_0,\sqrt{t})|}.$$
One concludes by using Assumptions $(A)$ and $(D)$. 

Let us now estimate $m(s, t, x_0)$.
First,  for $t \ge 1$,
\begin{eqnarray*}
c \sqrt{t}
\le r &= &\sqrt{c^{-1}_2 t (| \ln (C_1/s)|+1)}\\ & =&
\sqrt{c^{-1}_2 t (| \ln (4
C_1  F(t)/I(t))|+1)} \\
& \le&
\sqrt{c^{-1}_2 t (| \ln (
C  F(t)t^{\nu/2})|+1)}\\
&\le &C
\sqrt{t} \sqrt{\ln (e + t F(t))}.
\end{eqnarray*}
Finally,
\begin{equation}
m(s, t, x_0) \ge  \frac{c}{t
\ln (e+t F(t))} \left( \frac{I(t)|B(x_0,\sqrt{t})|}{  F^2(t) }
\right)^{2/\nu}. \label{clear}
\end{equation}

By \eqref{bec} and
\eqref{clear}, it follows that
\[
I'(t) \le -c  \frac{
I(t)^{1+(2/\nu)} |B(x_0,\sqrt{t})|^{2/\nu}}{t F(t)^{4/\nu} \ln(e+t
F(t))},
\]
that is
\begin{equation}
\frac{I'(t)}{I(t)^{1+(2/\nu)}}  \le -\frac{c
|B(x_0,\sqrt{t})|^{2/\nu}}{t F(t)^{4/\nu} \ln(e+t F(t))}.
\label{fro}
\end{equation}
Integrating (\ref{fro}) from $t/2$ to $t$ and
using the monotonicity of
$F(t)$ and $|B(x_0,\sqrt{t})|$, one easily
obtains
\begin{equation}
p^V_{2t}(x_0, x_0)=I(t) \le C \frac{ F^2(t) \left[\ln
(e+t
F(t))\right]^{\nu/2}}{|B(x_0,\sqrt{t})|}.
\label{just}
\end{equation}

From this on-diagonal bound, we can derive the full bound
\eqref{im}
by
either the
method in \cite{G2} or the wave equation method in \cite{Sik}.

\qed

\bigskip

\maketitle
\subsection{The $L^1$ to $L^1$
estimates}\label{1}

\begin{proposition}\label{1,1}
Suppose that  $M$
satisfies Assumptions (A), (B), (C),  (D), and that the
operator $\Delta -V$ is strongly positive.
If  $V \in L^p( M,\mu)$  for some $p\in[1,+\infty)$, then there exists $C=C(p)$ such that 
\[
\begin{cases}
\Vert
P_t^V \Vert_{1, 1} \le C t^{1/2}, \ \forall\, t \ge
1, \, \mbox{ if }p=1,\\
\Vert
P_t^V \Vert_{1, 1} \le C t^{(p-\theta)/2}, \ \forall\, t \ge
1, \, \mbox{ if }1<p<2,  \mbox{ for some } \theta=\theta(p)> 0,\\
\Vert
P_t^V \Vert_{1, 1} \le C t^{(p-1)/2}, \ \forall \,t \ge 1, \mbox{ if } p \ge
2.
\end{cases}
\]
\end{proposition}

{\it Remark} 
We shall see 
in Section \ref{wd} below that, if one weakens $(A)$, and $(B)$ as suggested in the remark after Proposition \ref{prel}, one can still prove
\[
\Vert P^V_t \Vert_{1, 1} \le  C_p t^{p/2}, \ \forall\, t \ge
1,\ 1\le p<+\infty .
\]
If  now one  only assumes that  $\Delta-V$ is positive instead of being strongly positive, 
then one can still prove, for  $V \in L^p( M,\mu)$ and $t\ge 1$,
\[
\Vert P^V_t \Vert_{1, 1} \le \begin{cases}  C t,
\mbox{ if } 1
\le p \le 2;\\
C t^{p/2},\mbox{ if }  p \ge 2;
\end{cases}
\]
Note that, according to  \cite{Si:1}, Theorem 3.1, the above estimate is sharp in the range 
$1\le p\le 2$, which shows the role of strong positivity
in the better estimate of Proposition \ref{1,1}.
For $p>2$,  similar estimates were proved in \cite{Li}, Theorem 8.1,
by extending the method of \cite{DS}, Theorem
3,  to the manifold
case.
This is also what we shall do to prove the more precise Proposition \ref{1,1},
taking advantage in addition of the strong positivity of $V$.

\medskip

\proof

Let us first prove that there
exists $\delta>0$ such that
\begin{equation}
\sup_y
\Vert p^V_t(., y) \Vert_2 = \Vert P^V_t \Vert_{2, \infty} \le C
t^{-\delta},
\quad t \ge 1. \label{2.9}
\end{equation}
The proof  goes as follows.

Given $f \in
C^\infty_0( M)$ and $q>1$, by  an easy consequence of the Feynman-Kac
formula (see for instance \cite{Z0}, p.712), one may write
\begin{equation}
|P_t^V f (x)|
\le [ e^{-t (\Delta - q V)} |f| (x)]^{1/q} [e^{t
\Delta} |f|(x)]^{(q-1)/q}.
\label{2.10}
\end{equation}
Since $\Delta -V$ is
strongly positive, when $q$ is sufficiently close to $1$, the
operator $\Delta - q V$ is also strongly positive, thus $P_t^{qV}$ is contractive on $L^2(M,\mu)$.
Hence
\begin{equation}
e^{-t (\Delta - q V)} |f| (x)
\le \|P_t^{qV}\|_{2,\infty}
\ \Vert f \Vert_2 \le \|P_{t-1}^{qV}\|_{2,2}\|P_1^{qV}\|_{2,\infty} \Vert f
\Vert_2\le\|P_1^{qV}\|_{2,\infty} \Vert f
\Vert_2,
\label{2.11}
\end{equation}
for $t\ge 1$.
Using the assumption that $p_t$ is
bounded from above by a
Gaussian (Assumption $(B)$), we have
\[
e^{-t \Delta}
|f|(x) \le \sup_y  \Vert p_t(y, .) \Vert_2 \ \Vert
f \Vert_2 \le
\frac{C}{\sqrt{|B(x,\sqrt{t})|}} \Vert f \Vert_2.
\]
By Assumptions $(A)$ and $(D)$,
\[
|B(x,\sqrt{t})|\ge
c_1 t^{\nu/2} |B(x,1)| \ge c_2 t^{\nu/2}.
\]
Therefore
\[
e^{t \Delta} |f|(x) \le
\frac{C}{t^{\nu/4}} \Vert f \Vert_2.
\]Substituting this and (\ref{2.11}) into
(\ref{2.10}), we have, for some
$\delta>0$,
\[
|P_t^V f (x)| \le
\frac{C}{t^{\delta}} \Vert f \Vert_2, \quad t \ge
1,\,x\in M.
\] This proves
(\ref{2.9}).

Now we are ready to prove the decay estimates in $L^1-L^1$ norm.

\bigskip

{\it
Case 1.} Assume $V
\in L^1(
M,\mu)$. 

Fix $y\in M$, and let $u(x,t)=p_t^V(x,y)$.
Since $u$ satisfies
\[
\Delta u - V u + u_t =0, 
\]
integrating on $[1,t]\times M$ yields
$$
\int_ M u(x, t) \,d\mu(x) = \int_ M u(x,
1) d\mu(x) +
\int^t_1\int_ M V(x) u(x, s)
\,d\mu(x)ds.
$$

By \eqref{small} and doubling (again, subexponential growth is enough, see the remark after Proposition \ref{prel}), it follows that
\begin{equation}
\int_ M u(x, t) \,d\mu(x)\le C+
\int^t_1\int_ M V(x) u(x, s)
\,d\mu(x)ds.\label{int}
\end{equation}

From \eqref{int} and the strong positivity of $V$,
\[
\aligned
 \int_ M u(x, t) d\mu(x)& \le  C+ \int^t_1\int_M
V(x) u(x, s) d\mu(x)ds\\
&\le C + \left(\int^t_1\int_ M V(x) d\mu(x)ds
\right)^{1/2}
\left(\int^t_1\int_ M V(x) u^2(x, s) d\mu(x)ds
\right)^{1/2}\\
&\le  C+ \sqrt{A} \sqrt{t} \Vert V
\Vert^{1/2}_1
\left(\int^t_1\int_ M |\nabla u|^2(x, s) d\mu(x)ds
\right)^{1/2}.
\endaligned
\]
Multiplying by $u$ the
equation
\[
\Delta u - V u + u_t =0
\]
and integrating on $[1,t]\times M$, we obtain
\[
\int^t_1\int_ M
|\nabla u|^2(x, s) d\mu(x)ds - \int^t_1\int_M V u^2 d\mu(x)ds +\frac{1}{2}
\int_ M u^2(x, t) d\mu(x) =\frac{1}{2}
\int_ M u^2(x, 1) d\mu(x).
\]Since
$\Delta -V$ is strongly positive, we obtain
\begin{equation}
\int^t_1\int_ M |\nabla
u|^2(x, s) d\mu(x)ds \le
\frac{1}{2(1-A)} \int_ M u^2(x, 1) d\mu(x)
\le
\frac{C}{2(1-A)}. \label{2.12}
\end{equation}
Finally
\[
\|P_t^V\|_{1\to 1}=\int_ M u(x, t)\, d\mu(x)
\le C \sqrt{\frac{A}{1-A}} \Vert V
\Vert^{1/2}_1 \sqrt{t},
\]
which is the claim.

\medskip

{\it Case
2.} Now we assume that $V \in L^p( M,\mu)$, $p \in (1,
2)$.

The decay
estimate in this range of $p$ seems to be new even in the
Euclidean case.

Using
H\"older's inequality repeatedly, one has
\begin{eqnarray*}
&& \int_ M u(x, t)
d\mu(x) \le C+ \int^t_1 \int_ M
V(x) u(x, s) d\mu(x)ds\\
&=&C+ \int^t_1 \int_M V(x)^{p/2} V(x)^{1-(p/2)} u(x, s) d\mu(x)ds\\
&\le& C + \left(\int^t_1
\int_ M V(x)^p d\mu(x)ds \right)^{1/2}
 \left(\int^t_1 \int_ M
V(x)^{2-p} u^2(x, s) d\mu(x)ds
 \right)^{1/2}\\
&=&C + t^{1/2} \Vert V
\Vert^{p/2}_p  \left(\int^t_1 \int_ M
V(x)^{2-p} u(x, s)^{2 (2-p)} u(x,
s)^{2p-2} d\mu(x)ds  \right)^{1/2}\\
&\le& C + t^{1/2} \Vert V \Vert^{p/2}_p
 \left(\int^t_1 \int_M[V(x)^{2-p} u(x, s)^{2 (2-p)}]^{1/(2-p)}
d\mu(y)ds
 \right)^{(2-p)/2}\\
&\qquad& \times
 \left(\int^t_1 \int_ M u(x,
s)^{(2p-2)/(p-1)} d\mu(x)ds  \right)^{(p-1)/2}\\
&=&C + t^{1/2} \Vert V
\Vert^{p/2}_p  \left(\int^t_1 \int_
MV(x) u^2(x, s) d\mu(y)ds
 \right)^{(2-p)/2}
 \left(\int^t_1 \int_ M u^2(x, s) d\mu(x)ds
\right)^{(p-1)/2}\\
&\le& C + A^{(2-p)/2} t^{1/2} \Vert V \Vert^{p/2}_p \left(
\int^t_1 \int_ M |\nabla u(x, s)| ^2 d\mu(y)ds  \right)^{(2-p)/2}
\left(\int^t_1 \int_ M u^2(x, s) d\mu(x)ds  \right)^{(p-1)/2}.
\end{eqnarray*}
Using (\ref{2.12}),
we obtain
$$\int_ M u(x, t)
d\mu(x) \le C + C't^{1/2}\left(\int^t_1 \int_ M u^2(x, s) d\mu(x)ds  \right)^{(p-1)/2}
$$

From (\ref{2.9}), we
know that
$$
\int^t_1 \int_ M u^2(x, s) d\mu(x)ds  \le Ct^{1-\delta},
$$
and using  (\ref{2.12}),  we deduce 

\[
\int_ M u(x, t) d\mu(x) \le C + C'
t^{1/2+(1-2\delta)(p-1)/2} =
C t^{(p-\theta)/2},
\]
\medskip
with $\theta=2\delta(p-1)$.

\bigskip

{\it
Case 3.} The only remaining case is when $V \in L^p(
M, \mu)$, $p \ge
2$.
\medskip

This case may be skipped in a first reading; indeed, if one is prepared to replace $(p-1)/2$
by $p/2$ in the estimate, the simpler proof in Section \ref{wd} will do.

The reader may guess that the claimed estimate follows from the idea in \cite{DS}, p.99,
where a similar bound for the Schr\"odinger heat kernel was proven
in the Euclidean case. However in \cite{DS}, the authors use the
boundedness of $\Delta^{-1/2}$ from some $L^p$ space to
another. But it is known that this property is false for most open
manifolds. Therefore, we have to work considerably harder. We will
show that the inverse square root of the Laplacian on $ M
\times \reel^3$ is bounded from the space $L^{p_1} \cap L^{p_2}$
to $L^2$ for some $p_1, p_2$. Then we will use the idea in \cite{DS} to
get an $L^1$ to $L^1$ bound for a version of our Schr\"odinger semigroup acting  on $ N \equiv  M
\times \reel^3$. After integrating over $\reel^3$, we will
reach the desired  $L^1$ to $L^1$ bound for the Schr\"odinger semigroup acting on $M$.

Points in $ N$ will be denoted by $\tilde{x}= (x, x')$ and
$\tilde{y}= (y, y')$, ..., where $x, y \in  M$ and $x', y'
\in \reel^3$. The distance function on $ N$ is denoted by
$d(\tilde{x}, \tilde{y}) = d(x, y) +|x'-y'|$, and the Riemannian measure on $N$ by $d\tilde{\mu}$. The Laplace-Beltrami
operator on $ N$ will be denoted by  $\tilde{\Delta}$. Denote by $\tilde{P}_t$ 
the heat semigroup
on $ N$  and let  $\tilde{P}_t^V=e^{-t(\tilde{\Delta} - V)}$, where $V$ is
the function on $M$ defined in (1.1). The kernel for
$\tilde{\Delta}^{-1/2}$ is
\begin{equation}
{\Delta}^{-1/2}(\tilde{x}, \tilde{y}) =
\int^{\infty}_0
t^{-1/2} \tilde{p}_t(\tilde{x}, \tilde{y}) dt.
\label{2.13}
\end{equation}

It is clear that $ N$ satisfies Assumptions $(A)$, $(B)$, $(C)$ and $(D)$
just like $M$ does. Moreover, if  $\Delta-V$ is strongly
positive on $M$, then $\tilde{\Delta} - V$ is strongly positive on $
N$. Now suppose we can prove that
\begin{equation}
\Vert \tilde{P}^V_t \Vert_{1, 1} \le C t^{(p-1)/2}. \label{2.14}
\end{equation}
Then, since
\begin{equation}
\tilde{p}^V_t(\tilde{x}, \tilde{y}) = p_t^V(x, y)
\frac{c}{t^{3/2}} e^{-|x'-y'|^2/4t}, \label{2.15}
\end{equation}
we
will deduce that
\[
\Vert P_t^V \Vert_{1, 1} = \Vert P_t ^V\Vert_{\infty, \infty} = \Vert
P_t ^V1 \Vert_ \infty= \Vert {\tilde P}^V_t 1 \Vert_ \infty =\Vert
\tilde{P}^V_t \Vert_{\infty, \infty} =\Vert \tilde{P}^V_t
\Vert_{1, 1} \le C t^{(p-1)/2},
\]thus finishing case 3. The rest of the section is
devoted to proving (\ref{2.14}). It will be
divided into several
steps.
\medskip

{\it Step 1.} We show that there exists $C>0$ such
that
\begin{equation}
\tilde{\Delta}^{-1/2}(\tilde{x}, \tilde{y}) \le
C
\frac{d(\tilde{x}, \tilde{y})}{|B(\tilde{x}, d(\tilde{x},
\tilde{y}))|}.
\label{2.16}
\end{equation}
F{}rom the upper bound for $\tilde{P}_t$, which is a
consequence
of the Gaussian upper bound for $P_t$, we have
\[
\aligned I & \equiv
\tilde{\Delta}^{-1/2}(\tilde{x},
\tilde{y})\\
&\le C \int^\infty_0
\frac{1}{\sqrt{t}} \frac{ e^{-c
d^2(\tilde{x},
\tilde{y})/t}}{|B(\tilde{x},
\sqrt{t})|} dt \\
&=C
\int^{d^2(\tilde{x}, \tilde{y})}_0 \frac{1}{\sqrt{t}}
\frac{ e^{-c
d^2(\tilde{x}, \tilde{y})/t}}{|B(\tilde{x},
\sqrt{t})|} dt
+C \int^{\infty}_{d^2(\tilde{x},
\tilde{y})} \frac{1}{\sqrt{t}} \frac{
e^{-c d^2(\tilde{x},
\tilde{y})/t}}{|B(\tilde{x}, \sqrt{t})|}
dt\\
&\equiv C I_1 + C I_2.
\endaligned
\label{2.17}
\]By
the doubling condition, for $t \le
d^2(\tilde{x},
\tilde{y})$,
$$
|B(\tilde{x}, \sqrt{t})| \ge c |B(\tilde{x},
d(\tilde{x},
\tilde{y}))| \  \left(\frac{\sqrt{t}}{d(\tilde{x},
\tilde{y})}
\right)^\nu.
$$
 Hence 
\begin{equation}
I_1 \le C
\frac{d(\tilde{x}, \tilde{y})}{|B(\tilde{x},
d(\tilde{x}, \tilde{y}))|}
\label{2.18}
\end{equation}

Next we estimate $I_2$.

By (\ref{2.15}), we have
\[
I_2  \le C
\int^\infty_{d^2(\tilde{x}, \tilde{y})}
\frac{1}{\sqrt{t}} \frac{1}{|B(x,
\sqrt{t})| \ t^{3/2}} dt \le
\frac{C}{|B(x, d(\tilde{x},
\tilde{y}))|}
\int^\infty_{d^2(\tilde{x}, \tilde{y})}
\frac{dt}{t^2}
=\frac{C}{|B(x, d(\tilde{x}, \tilde{y}))| \
d^2(\tilde{x},
\tilde{y})}.
\]Note that
\[
|B(\tilde{x}, d(\tilde{x},
\tilde{y}))| = |B(x, d(\tilde{x},
\tilde{y}))| \ d^3(\tilde{x},
\tilde{y}).
\]The above implies
\begin{equation}
I_2 \le \frac{C d(\tilde{x},
\tilde{y})}{|B(\tilde{x},
d(\tilde{x}, \tilde{y}))|}. \label{2.19}
\end{equation}
The
combination of (\ref{2.19}) and (\ref{2.18})  implies (\ref{2.16}).
\medskip

{\it Step
2.} We prove that there exist $p_1, p_2>1$ and $C>0$
such that
\begin{equation}
\Vert
\tilde \Delta ^{-1/2} f \Vert_2 \le C ( \Vert f
\Vert_{p_1} + \Vert f
\Vert_{p_2}), \label{2.20}
\end{equation}
for all $f \in L^{p_1} \cap L^{p_2}.$

From
(\ref{2.16}),
\[
\aligned
 |\tilde \Delta ^{-1/2} f( \tilde{x})| &\le
C
\int_{\tilde M} \frac{ d(\tilde{x},
\tilde{y})}{|B(\tilde{x},
d(\tilde{x}, \tilde{y}))|} | f(
\tilde{y})| \,d\tilde{\mu}(\tilde{y})\\
&\le C \int_{d(\tilde{x}, \tilde{y}) \le 1}... +
C
\int_{d(\tilde{x}, \tilde{y}) \ge 1} \equiv C J_1(\tilde{x}) +
C
J_2(\tilde{x}).
\endaligned
\label{2.21}
\]

By Young's
inequality
\[
\Vert J_1 \Vert_2 \le C  \Vert f \Vert_{p_1}
\sup_{\tilde{x}}
\bigg{(} \int_{d(\tilde{x}, \tilde{y}) \le 1} \bigg{[}
\frac{
d(\tilde{x}, \tilde{y})}{|B(\tilde{x}, d(\tilde{x},
\tilde{y}))|}
\bigg{]}^{r_1} \,d\tilde{\mu}(\tilde{y}) \bigg{)}^{1/r_1}
\]where $(1/p_1) + (1/r_1) = 1 + (1/2)$. Hence
\[
\aligned \Vert  J_1
\Vert_2 &\le  \Vert f \Vert_{p_1}
\sup_{\tilde{x}} \ \sum^\infty_{k=0}
\bigg{(} \int_{2^{-(k+1)} \le
d(\tilde{x}, \tilde{y}) \le 2^{-k}} \bigg{[}
\frac{d(\tilde{x},
\tilde{y})}{|B(\tilde{x}, d(\tilde{x},
\tilde{y}))|}
\bigg{]}^{r_1} \,d\tilde{\mu}(\tilde{y})
\bigg{)}^{1/r_1}\\
&\le C \Vert f
\Vert_{p_1} \sup_{\tilde{x}} \ \sum^\infty_{k=0}
\frac{2^{-k
r_1}}{|B(\tilde{x}, 2^{-(k+1)})|^{r_1-1}}.
\endaligned
\]Using the doubling
property and the fact that, thanks to Assumption $(D)$,
$|B(\tilde{x}, 1)| \ge c>0$, we have
\[
|B(\tilde{x}, 2^{-(k+1)})| \ge C 2^{-k \nu} |B(\tilde{x},
1)| \ge
C 2^{-k \nu}.
\]Hence
\[
\Vert J_1 \Vert_2 \le  \Vert f \Vert_{p_1}
\sum^\infty_{k=0} 2^{-k
r_1} 2^{k \mu (r_1-1)}.
\]Choosing $r_1$
sufficiently close to $1$, the above series is
convergent.
Therefore
\begin{equation}
\Vert J_1 \Vert_2 \le  C \Vert f \Vert_{p_1}.
\label{2.22}
\end{equation}
Similarly, by Young's inequality again,
\[
\Vert J_2 \Vert_2 \le C  \Vert f \Vert_{p_2}
\sup_{\tilde{x}}
\bigg{(} \int_{d(\tilde{x}, \tilde{y}) \ge 1} \bigg{[}
\frac{
d(\tilde{x}, \tilde{y})}{|B(\tilde{x}, d(\tilde{x},
\tilde{y}))|}
\bigg{]}^{r_2} \,d\tilde{\mu}(\tilde{y}) \bigg{)}^{1/r_2}
\]where $(1/p_2) + (1/r_2) = 1 + (1/2)$. Hence
\[
\aligned \Vert J_2
\Vert_2 &\le C \Vert f \Vert_{p_2}
\sup_{\tilde{x}} \ \sum^\infty_{k=0}
\bigg{(} \int_{2^{k} \le
d(\tilde{x}, \tilde{y}) \le 2^{k+1}} \bigg{[}
\frac{ d(\tilde{x},
\tilde{y})}{|B(\tilde{x}, d(\tilde{x},
\tilde{y}))|}
\bigg{]}^{r_2} \,d\tilde{\mu}(\tilde{y})
\bigg{)}^{1/r_2}\\
&\le C \Vert f
\Vert_{p_2} \sup_{\tilde{x}} \ \sum^\infty_{k=0}
\frac{2^{k
r_2}}{|B(\tilde{x}, 2^{k})|^{r_2-1}}.
\endaligned
\]Using the fact that
$|B(\tilde{x}, 2^{k})| \ge c 2^{3k}$ which follows from
the definition of
$ N$ and Assumption $(D)$, we have
\[
\Vert J_2
\Vert_2 \le  \Vert f \Vert_{p_2} \sum^\infty_{k=0} 2^{k
r_2 - 3 k
(r_2-1)}.
\]Choosing $r_2$ sufficiently large, the above series
is
convergent. Therefore
\begin{equation}
\Vert J_2 \Vert_2 \le  C \Vert f \Vert_{p_2}.
\label{2.23}
\end{equation}
Inequality (\ref{2.20}) immediately follows  from (\ref{2.23})
and (\ref{2.22}).

\medskip
{\it Step 3.} As  in \cite{DS},  by Duhamel's formula, one has
\begin{eqnarray*}
 \Vert \tilde{P}^V_{t+1} \Vert_{\infty, \infty} &=&
\Vert \tilde{P}^V_{t+1} 1 \Vert_{\infty}\\
& \le  &\Vert \tilde{P}^V_{1} 1 \Vert_{\infty}+
\int^t_1 \Vert \tilde{P}^V_{s+1} V
\Vert_{\infty} \,
ds\\
  & \le  &C+
\int^t_1 \Vert \tilde{P}^V_{s+1} V
\Vert_{\infty} \,ds,
\end{eqnarray*}
hence by    interpolation, 
\begin{equation} \Vert \tilde{P}^V_{t+1} \Vert_{\infty, \infty} \le C + \int^t_1 \Vert \tilde{P}^V_{s+1} \Vert^{2/p}_{2, \infty}
\Vert \tilde{P}^V_{s+1} \Vert^{1-(2/p)}_{\infty, \infty} \Vert V
\Vert_{p}\, ds.
\label{2.24}
\end{equation}

Let us now estimate $\Vert \tilde{P}^V_{s+1} \Vert^{2/p}_{2, \infty}$.
 Using the strong positivity of $\tilde{\Delta}-V$  on $N$, one has
\[
(\tilde{\Delta}
-V) \ge a^2 \tilde{\Delta}
\]for some $a>0$. Therefore, for all $f \in
C^\infty_0( N)$,
\[
\Vert (\tilde{\Delta} -V)^{-1/2} f \Vert_2 \le
a^{-1} \Vert
{\tilde{\Delta}}^{-1/2} f \Vert_2 \le C \Vert f
\Vert_X.
\]Here $X = L^{p_1} \cap L^{p_2}$ where $p_1, p_2$ are given
in
(\ref{2.20}) and 

\begin{equation}
\Vert f \Vert_X = \Vert f \Vert_{p_1} + \Vert
f
\Vert_{p_2}.\label{x}
\end{equation}

 Hence
\[
\aligned \Vert e^{- t (\tilde{\Delta} -V) } f
\Vert_2 &= t^{-1/2}
\Vert e^{-(\tilde{\Delta} -V)
t} ((\tilde{\Delta} -V)
t)^{1/2} (\tilde{\Delta} -V)^{-1/2} f \Vert_2 \\
&\le C t^{-1/2} \Vert
(\tilde{\Delta} -V)^{-1/2} f \Vert_2 \le
\Vert
\tilde{\Delta}^{-1/2} f
\Vert_2\\
&\le C  t^{-1/2} \Vert f \Vert_X. \endaligned
\]This
shows
\[
\Vert \tilde{P}_{t} \Vert_{2, X^*} = \Vert \tilde{P}_{t} \Vert_{X, 2} \le C t^{-1/2},
\]where $X^*$ is the
dual of $X$. Then write
\begin{equation}
\Vert \tilde{P}_{s+1}^V \Vert_{2, \infty} \le  \Vert
\tilde{P}_{s}^V
\Vert_{2, X^*} \Vert \tilde{P}_{1} ^V\Vert_{X^*,
\infty}.
\label{2.25}
\end{equation}
It follows easily from the bound
\[
\tilde{p}_1^V(\tilde{x}, \tilde{y}) \le \frac{C}{|B(\tilde{y},
1)|} e^{-c d^2(\tilde{x}, \tilde{y})},
\]
and Assumptions $(A)$ and $(D)$
 that
 $\tilde{P}_1^V$ is bounded from $L^1$ to any $L^p$, $1\le p\le +\infty$.
 Therefore
 \begin{equation}
\Vert \tilde{P}_{1}^V \Vert_{X^*, \infty}=\Vert \tilde{P}_{1}^V \Vert_{1, X} =\Vert \tilde{P}_{1}^V \Vert_{1, p_1}+\Vert \tilde{P}_{1}^V \Vert_{1, p_2}<+\infty.
\label{2.26}
\end{equation}
This  together with
(\ref{2.25})
implies
\[
\Vert \tilde{P}^V_{s+1} \Vert_{2, \infty} \le C
s^{-1/2}.
\]
 Using this and
(\ref{2.24}), we obtain
\[
\Vert \tilde{P}^V_{t+1} \Vert_{\infty, \infty} \le C + C \Vert V
\Vert_p \int^t_1 s^{-1/p}
  \Vert \tilde{P}^V_{s+1}
\Vert^{1-(2/p)}_{\infty, \infty} ds.
\]

By Gronwall's lemma, 
\[
 \Vert \tilde{P}^V_{t+1} \Vert_{\infty, \infty}  \le C t^{p/2-1/2}.
\]
Since\[
\Vert \tilde{P}^V_{t+1} \Vert_{1, 1} = \Vert \tilde{P}^V_{t+1} \Vert_{\infty, \infty},
\]
 the proof
of  Proposition \ref{1,1} is complete. \qed

\maketitle
\subsection{Proof of Theorem \ref{mainV}}\label{end}

Under the assumptions of Theorem \ref{mainV}, 
the combination of Propositions \ref{nash} and \ref{1,1} yields
\begin{equation}
p^V_t(x, y)\le  \frac{Ct \ln^{\nu/2}(e+ t)}{|B(x,\sqrt{t})|} \exp (- c d^2(x, y)/t), \ \forall\,x,y\in M, \, t\ge 1
\label{step}
\end{equation}
if $p=1$,
\begin{equation}
p^V_t(x, y)\le  \frac{Ct ^{p-\theta}}{|B(x,\sqrt{t})|} \exp (- c d^2(x, y)/t), \ \forall\,x,y\in M, \, t\ge 1
\label{step1}
\end{equation}
for some $\theta=\theta(p)>0$, if $1\le p<2$, and
\begin{equation}
p^V_t(x, y)\le  \frac{Ct ^{p-1}\ln^{\nu/2} (e+t)}{|B(x,\sqrt{t})|} \exp (- c d^2(x, y)/t), \ \forall\,x,y\in M, \, t\ge 1
\label{step2}
\end{equation}
if $p\ge 2$.

The above bound does not  explicitely reflect the contribution of
the constant $A$
(which measures the size of $V$).  To remedy this, we use
again \eqref{2.10}, which yields
\begin{equation}
p_t^V(x, y) \le \left( p^{V'}_t(x, y)\right)^{1/q} \left(p_t(x,
y)
\right)^{(q-1)/q}, \forall\,x,y\in M, \,t>0,\label{fek}
\end{equation}
where $V'=qV$ and $q>1$.

Assume $A>0$, otherwise there is nothing to prove. If $q<1/A$, then
$V'$ is obviously strongly positive with constant $qA<1$.

If $1<p<2$,  \eqref{step1}  yields
\begin{equation}
p^{V'}_t(x, y) \le 
\frac{Ct^{p-\theta}}{|B(x,\sqrt{t})|} \exp
(- c d^2(x, y)/t).
\label{v1}
\end{equation}

Applying (\ref{v1})
and using the Gaussian upper bound on $p_t$, we
obtain
$$
p_t^V(x, y) \le  \frac{ C t^{(p-\theta)/q}
}{|B(x,\sqrt{t})|} \exp (- c d^2(x, y)/t).
$$Taking $q$ sufficiently close
to
$1/A$,  one can make $(p-\theta)/q<pA$, which finishes the proof of
  Theorem \ref{mainV} in the case $1<p<2$. The proofs when $p=1$ and $p \ge
 2$ are identical using \eqref{step}  and \eqref{step2}  instead of \eqref{step1}.
\qed
\medskip

\maketitle
\section{Gaussian bound on the heat kernel on forms and  the Riesz transform}\label{rie}
\medskip

Our next theorem provides an all time Gaussian upper bound for the heat kernel on
$1$-forms on complete non-compact Riemannian manifolds satisfying Assumptions
$(A)$, $(B)$, $(C)$, together with a certain condition of smallness of the Ricci curvature.
Using this  theorem and an argument in \cite{CD}, one deduces a proper
bound
for the gradient of the heat kernel on functions. By the main result in
\cite{ACDH}, one obtains the $L^p$ boundedness of the Riesz transform on
these
manifolds for all $1<p<+\infty$.  Let us point out that in \cite{Li}, Theorem 9.1, another sufficient condition in terms of Ricci curvature is given for    $L^p$ boundedness of  the Riesz transform. However this condition seems to exclude Ricci curvature bounded from
below together with non-compactness.

\begin{theorem}\label{riet}
Let ${M}$ be a complete non-compact Riemannian manifold satisfying Assumptions (A), (B) and (C).  Then there exists $\delta>0$ depending
only
on the constants in (A) and (B) such that
for
any $c > 1/4$, there exists  $C>0$ such that
\[
|\vec{p}_t(x, y)| \le \frac{C}{|B(x, \sqrt{t})|}\exp\left( - c
d(x,
y)^2/t\right),
\]for all $x, y \in { M }$ and $t >0$, provided that
\begin{equation}
K(V) \equiv \sup_{x \in  { M }} \int^\infty_0 \int_{  M }
\frac{1}{|B(x, \sqrt{s})|} e^{ -  d^2(x, y)/s} V(y) d \mu(y) ds <
\delta.\label{coco}
\end{equation}
\end{theorem}

As we already said, the following statement is a consequence from Theorem \ref{riet} and either Theorem 5.5 in \cite{CD}, or  Theorem 1.4
in \cite{ACDH} together with \cite{CD}, pp.1740-1741.
This extends the class of manifolds for which one can answer a question asked by Strichartz in \cite{St}.

\begin{corollary}
Let ${M}$ be a complete non-compact Riemannian manifold satisfying Assumptions (A), (B), (C),
and condition \eqref{coco} for $\delta>0$ small enough. 
Then, for all $p\in (0,+\infty)$, there exist $C_p,c_p>0$ such that
$$c_p\||\nabla f|\|_p\le \|\Delta^{1/2} f\|_p\le C_p\||\nabla f|\|_p,\ \forall\,f\in C^{\infty}_0( M). $$
\end{corollary}
\noindent

{\it Proof of Theorem \ref{riet}.}
As we have seen in the introduction,
\[
|\vec{p}_t(x, y)| \le p^{V}_t(x, y)
\]under Assumptions $(A)$ to $(C)$.  Now let us recall Theorem A, part (b)
in
\cite{Z}.  It implies
that
\[
  p^{V}_t(x, y) \le
\frac{C}{|B(x, \sqrt{t})|} e^{ - c d^2(x, y)/t},
\]for all $x, y \in {M }$ and $t >0$, provided that, for certain
$c_0,
\varepsilon_0>0$, there holds
\[
N(V)  < \varepsilon_0.
\]Here
\begin{eqnarray*}
N(V) &\equiv &
  \sup_{x \in  { M }, t>0} \int^t_{0} \int_{ M }
\frac{e^{ - c_0 d^2(x, y)/(t-s)} }{|B(x, \sqrt{t-s})|} V(y) d 
\mu(y)
ds  \\
&&\qquad  + \sup_{y \in  { M }, s>0} \int^{\infty} _s \int_{ M }
\frac{e^{ - c_0 d^2(x, y)/(t-s)} }{|B(x, \sqrt{t-s})|} V(x) d 
\mu(x)
dt.
\end{eqnarray*}

We should mention that this theorem was stated for a doubling metric in
the Euclidean space and for time dependent functions $V$, under the
extra
assumption $(D)$. However the proof was a general one applicable verbatim to any
manifold under Assumptions $(A)$ and $(B)$ only.

Changing  variables and using doubling, one sees that
\[
N(V)  \le C_0  K(V),
\] 
where $C_0$ only depends on the doubling constants.
The conclusion follows.
\qed

{\it Remarks:}

- Suppose in addition that $|B(x, r)| \sim r^n$ with $n>2$, uniformly in $x\in M$, for  large $r$, then it is any easy exercise to check that the theorem
holds if
\[
V(x) \le \frac{a}{1 + d(x,x_0)^{2 + b}}
\]with any $b>0$ and $a$ sufficiently small, for some fixed $x_0\in M$.

- Let $({  M}, g_0)$ be a nonparabolic manifold with
nonnegative
Ricci curvature and volume growth property as in the last remark. Let
$h$
be another metric and $\eta$ be a smooth cut-off function on ${  
M}$.
Then the manifold   $({  M}, g)$ with $g=g_0 + \lambda \eta h$  is
covered by the theorem when $\lambda\ge 0$ is sufficiently small.
This is so because the constants in $(A)$ and $(B)$ are uniformly bounded
when $ 0 \le \lambda \le 1$ while $V=V(x)$ (for the metric $g$),  being
a
compactly supported function is arbitrarily small when $\lambda \to 0$.
\bigskip

\maketitle
\section{The case of non-negative  
Ricci curvature outside a compact set}\label{bad}
\medskip

In the next theorem, we establish an upper bound for the heat kernel on 
$1$-forms assuming Ricci curvature is nonnegative outside a compact set. The
upshot of the theorem
is that no other restriction on the Ricci curvature is needed.
This
upper bound
gives a good control of the heat kernel even in the presence of 
harmonic
forms.
In general one can not expect the heat kernel on forms to decay to zero, due to the possible presence of $L^2$ harmonic forms. Here we are able to 
show that the heat kernel has certain spatial  decay anyway. Using the spectral decomposition
of heat kernels, one can see that the upper bound in the theorem below is quite sharp near 
the diagonal at least. As far as we know, this bound is new even for Schr\"odinger 
heat kernels in the Euclidean case.

Moreover the assumption that the Ricci curvature is $0$ outside of a compact set 
can be improved to assuming that the negative part of the Ricci curvature decays
sufficiently fast near infinity. But we will not seek the full generality this time.

\begin{theorem}
Let ${M}$ be a manifold satisfying Assumptions (A), (B),  (C) and
(D).
In addition, we assume that the Ricci
curvature of ${  M}$ is nonnegative outside a compact set and the
manifold is
nonparabolic.
Then, for a fixed $0\in M$, there exist  $C,c>0$ such that
\[
|\vec{p}_t(x, y)| \le C  \min \{ \Gamma(x, 0) \Gamma(y, 0), 1
\}
e^{ - c
d^2(x, y)/t} +  \frac{C}{ |B(x, \sqrt{ t})| }
e^{ - c d^2(x, y)/t}
\]for all $x, y \in {  M }$, $t >1$.
Here $\Gamma$ is the Green's function of the  Laplacian $\Delta$ on ${ 
M}$.
\proof
\end{theorem}

We divide the proof into two steps.

{\it Step 1.}

As in Section \ref{poin}, we need a preliminary estimate.
\begin{equation}
| \vec{p}_t(x, y) | \le C e^{-c d^2(x, y)/t},
\ t\ge
1,
\label{expv}
\end{equation}
for some $C,c>0$.
Note that this estimate does not follow from Proposition \ref{prel} and \eqref{dom},
since we do not assume strong positivity of $\Delta-V$ any more.

However, the method is very similar to the one in Proposition  \ref{prel}. Let $u_0$ be a smooth compactly supported $1$-form. Write
\[
u(x, t) = \vec{P}_t
u_0(x).
\]Direct computations show, for any fixed $y \in  M$ and
$D>0$,
\[
\aligned
   \frac{d}{dt} &\int_ M |u|^2e^{ \frac{d^2(x,
y)}{D t}} \,d\mu(x)\\
   &=-2 \int_ M e^{\frac{d^2(x, y)}{D t}}  u\cdot \vec{\Delta} u
\,d\mu(x)
   - \int_ M |u|^2 e^{ \frac{d^2(x, y)}{D t} } \frac{d^2(x,
y)}{D
   t^2} \,d\mu(x).
\endaligned
\]Noticing that $\vec{\Delta}= d^* d + d d^*$,
the above implies, after
integration by parts,
\[
\aligned
   \frac{d}{dt}
&\int_ M |u|^2e^{ \frac{d^2(x, y)}{D t}} \,d\mu(x)\\
   &= - 2 \int_M  e^{\frac{d^2(x, y)}{D t}}\left(d \left(\frac{d^2(x, y)}{D 
t}
   \right)
\wedge
     u\right)\cdot d u  \,d\mu(x)  -2  \int_ M  e^{\frac{d^2(x, y)}{D t}}
du\cdot
     du \,d\mu(x)  \\
     & \qquad -2 \int_ M  e^{\frac{d^2(x, y)}{D
t}}| d^*u|^2
    \,d\mu(x) - \int_ M |u|^2 e^{ \frac{d^2(x, y)}{D
t} }
\frac{d^2(x, y)}{D
   t^2} \,d\mu(x)\\
   &\le C \int_ M
e^{\frac{d^2(x, y)}{D t}}\frac{d(x, y)}{D t}
   |u| |du| \,d\mu(x)  -
2 \int_ M  e^{\frac{d^2(x, y)}{D t}} |du|^2 \,d\mu(x)\\
&\qquad
     -  \int_M |u|^2 e^{ \frac{d^2(x, y)}{D t} } \frac{d^2(x,
y)}{D
   t^2} \,d\mu(x).
\endaligned
\]Using the inequality $\frac{d(x,
y)}{D t} |u| |du| \le
\frac{1}{\varepsilon} \frac{d^2(x, y)}{D^2 t^2}  |u|^2 +  \varepsilon
|du|^2$, we find
that
\[
\frac{d}{dt} \int_ M |u|^2 e^{ \frac{d^2(x,
y)}{D t}} \,d\mu(x) \le
0
\]when $D$ is sufficiently large.

Letting  $u_0$ converge to the
Dirac delta function
centered at $y$, we obtain
\[
\int_ M | \vec{p}_t(x, y)|^2   e^{
\frac{d^2(x, y)}{D t}} \,d\mu(x)
\le \int_ M | \vec{p}_1(x, y)|^2   e^{
\frac{d^2(x, y)}{D}} \,d\mu(x)
\]when $t \ge 1$. By the semigroup
domination property \eqref{dom}, 
Assumptions $(C)$ and $(B)$,
\[
| \vec{p}_1(x,
y)| \le  p_1^V(x, y) \le C p_1(x, y)  \le
\frac{C}{|B(x,1)|} e^{-c
d^2(x, y)}.
\]Integrating and using Assumptions $(D)$ and $(A)$, we have, for a suitable $D>0$,
\begin{equation}
\int_ M |
\vec{p}_t(x, y)|^2   e^{ \frac{d^2(x, y)}{D t}} \,d\mu(x)
\le C, \quad t \ge
1.\label{coll}
\end{equation}

Next, using the semigroup
property
\[
\aligned |\vec{p}_{2t}(x, y)| &= |\int_ M \vec{p}_t(x,
z)
\vec{p}_t(z, y) d\mu(z)|
\\
&=|\int_ M e^{ \frac{d^2(x, z)}{2 D t}}
\vec{p}_t(x, z) e^{
\frac{d^2(z, y)}{2D t}} \vec{p}_t(z, y) e^{-\frac{d^2(x,
z)}{2 D
t} -
\frac{d^2(z, y)}{2 D t}} d\mu(z) |\\
&\le  e^{ -\frac{d^2(x, y)}{4
D t}} \big{[} \int_ M e^{
\frac{d^2(x, z)}{ D t}} |\vec{p}_t(x, z)|^2
d\mu(z) \big{]}^{1/2}
\big{[} \int_ M e^{ \frac{d^2(y, z)}{ D t}}
|\vec{p}_t(y,
z)|^2 d\mu(z) \big{]}^{1/2}.
\endaligned
\]
Together with
\eqref{coll}, this  implies
\eqref{expv}.
\medskip

{\it Step 2.}

We assume that the Ricci curvature is nonnegative outside of a ball
$B(0,
A)$ for a fixed $A>0$.
Write, for any given $x_0 \in \partial B(0, A)$,
\[
u(y, t) = |\vec{p}_t( x_0,y)| .
\]Then it is an immediate consequence of  Bochner's formula (see for instance \cite{DL}, Lemma 4.1) that $u$ is a subsolution of the scalar heat equation in $B^c(0, A) 
\times
(0, \infty)$. i.e.
\[
\Delta u(y, t) +u_t(y, t) \le 0 .
\]Since, according to \eqref{expv}, $u$ is bounded from above by a constant for $t\ge 1$, and 
$\Gamma(y, 0)$ is bounded from below by a positive constant on any compact set,
there exists $C>0$ such that
\[
u(y, t) \le C   \Gamma(y, 0),\ \forall\,  y \in \partial
B(0,
A), \,t\ge 1.
\]Moreover, using again \eqref{expv},  
$$
u(y, 1) \le C e^{- c d^2(y, x_0)} \le C' \Gamma(y, 0)
$$
for
$y \in
  B^c(0, A)$.
The last inequality is due to the Cheng-Yau  gradient estimate from \cite{CY}, which
implies 
$$
\Gamma(y, 0)
\ge c e^{- C d(y,0)}
$$for some positive constants $C,c>0$
and 
$y \in
  B^c(0, A)$.
Now by the maximum principle, using  \eqref{expv} again, we deduce
\[
u(y, t)   \le C  \Gamma(y, 0)
\]for all $t \ge 1$ and $y \in B^c(0, A) $.  Here we just used the simple
observation that $ \Gamma(y, 0)$
is a solution of the scalar heat equation, whereas $u$ is a subsolution as already observed.

We have proved that
\[
|\vec{p}_t( x_0,y)|  \le C  \Gamma(y, 0)
\]for all $t \ge 1$, $y \in B^c(0, A) $ and $x_0 \in \partial B(0, A)$.
Let us now to explain how to keep the same estimate while moving away $x_0$.

For a fixed $y \in B^c(0, A)$, define the function
\[
w(x, t) = |\vec{p}_t(x, y)| .
\]Then $w$ is a subsolution of the scalar heat equation on $B^c(0, A)\times (0,+\infty)$.
For $x  \in \partial B(0, A)$, by the above estimate on $u$ we have
\[
w(x, t) = |\vec{p}_t(x, y)|  \le C  \Gamma(y, 0).
\]
Since $\Gamma(x, 0)$ is bounded away from $0$ for  $x  \in \partial B(0, A)$,  it holds
\[
w(x, t) \le C' \Gamma(y, 0)  \Gamma(x, 0)
\]
for some $C'>0$. It is clear that the function
\[
h(x, t) = \int_M p_{t-1}(x, z) w(z, 1) d\mu(z) +  C  \Gamma(y, 0)  \Gamma(x, 0)
\]is a solution of the scalar heat equation in $ B^c(0, A)  \times [1,
\infty)$.   Moreover, on the parabolic
boundary of the region, $h$ dominates $w$. By the maximum principle
again
\[
w(x, t) \le h(x, t) = \int_M p_{t-1}(x, z) w(z, 1) d\mu(z) +  C  \Gamma(y, 0)
\Gamma(x, 0)
\]for $x \in B^c(0, A)$ and $t \ge 1.$
Next we estimate the above  integral
term in the following way, by using the Gaussian upper bound for $p_t$:
\begin{eqnarray*}
\int_M p_{t-1}(x, z) w(z, 1) d\mu(z)  &\le& \int p_{t-1}(x, z) e^{- c d(y,  z)^2}
d\mu(z)\\
&\le&  \int_{d(x, z) \ge d(x, y)/2} ... d \mu(z)   +  \int_{d(y, z) \ge d(x,
y)/2} ...  d\mu(z) \le
  \frac{C}{ |B(x, \sqrt{ t})| } .
\end{eqnarray*}

Finally, incorporating \eqref{expv}, 
\[
 |\vec{p}_t(x, y)| = w(x, t) \le C \min \{ \Gamma(x, 0) \Gamma(y, 0), 1
\}
 +  \frac{C}{ |B(x, \sqrt{ t})| }
\]for all $x, y $ in $M$ and $t \ge 2$. This is the desired on-diagonal estimate.

Now the theorem follows  from the standard process of going from on to off-diagonal estimate. See \cite{Sik} e.g.\qed

\maketitle
\section{Bounds on manifolds without doubling condition}\label{wd}

In this final section we turn to noncompact manifolds not necessarily satisfying 
the volume doubling condition. This class of manifolds offers a much richer variety
than the doubling ones.

We show that under reasonable conditions the on-diagonal upper bound on the heat kernel on forms differs from that on functions only by a suitable power of time $t$.

 Let us consider a $n$-dimensional  manifold $M$ with Ricci curvature bounded from below, and whose small balls do not collapse, in other words Assumptions $(C)$ and $(D)$ are satisfied.
 Then 
\begin{equation}
p_t(x,y)\le Ct^{-n/2}\exp\left(-cd^2(x,y)/t\right),\ \forall\,  0<t\le 1,\,x,y\in M,\label{smol}
\end{equation}
for some $C,c>0$ and $p_t^V(x,y)$, $|\vec{p}_t(x, y)| $ satisfy similar estimates.

\medskip

Let us assume that the heat kernel on functions has a uniform rate of decay $\gamma$, where $\gamma$ is increasing, $C^1$ and one-to-one on $\reel_+$:
\begin{equation}
\sup_{x\in M}p_t(x, x) \le  \frac{
1}{\gamma( t)}, \ \forall\,\ t >0.
\label{gam}
\end{equation}
According to \cite{G}, this implies the following so-called uniform Faber-Krahn inequality:
 for any set $\Omega \subset  M$,
\[
\lambda_1( \Omega) \ge
\Lambda ( |\Omega|),\leqno{(UFK)}
\]
where $\Lambda$ is given by
\[
\Lambda(t)=\frac{\gamma'(t)}{\gamma(t)},
\]
i.e.
\[
t = \int^{\gamma(t)}_0 \frac{
d\eta}{\eta \Lambda(\eta)}.
\]
Conversely, if $\gamma$ satisfies a mild condition, the converse is true.
For more on this as well as examples where one can compute $\Lambda$, therefore $\gamma$, see for instance \cite{Cs}.

\medskip

We can now state our result in this setting.

\begin{theorem}\label{final}
Suppose $M$ satisfies  Assumptions (C) and (D). Assume that $V \in
L^p(M, \mu)$ for some $p\in[1,+\infty)$, and
that  $\Delta - V$ is strongly positive.  Finally assume that the heat kernel on functions
on $M$ satisfies the estimate \eqref{gam}. Then there exist positive
constants $c$ and $C$ such that
\[
|\vec{p}_t(x, x)| \le  \frac{C
t^{p}}{\gamma(c t)}, \ \forall\, \ t \ge 1, \,x\in M.
\]
\end{theorem}

\proof
We divide the proof into two parts.
\medskip

{\it Step 1. } $L^1$ to $L^1$ bound. 

Under the assumptions of the theorem we will prove that
\begin{equation}
\Vert
P_t^V \Vert_{1, 1} \le C t^{p/2}, \ \forall \,t \ge 1, \mbox{ if } p \ge
1.\label{1to1}
\end{equation}

Comparing with the proof of Proposition \ref{1,1}, we no longer 
have the doubling condition. However the growth rate of $\Vert
P_t^V \Vert_{1, 1}$ here is worse. On the other hand,  the proof, a simple application of the idea of \cite{DS}, is much shorter. 
In case $p>2$, the above estimate is essentially contained in \cite{DS}, see also \cite{Li}. We present the proof for completeness.
\medskip

{\it Case 1.} Assume $V \in L^1(M)$.

Fixing $y$, we write
$$
u(x, t) =  p^V_t(x, y).
$$

In this case the proof is almost identical to that of Case 1, Proposition \ref{1,1}.
The only change is that we use the small time bound \eqref{smol} on $p_t^V$ and the subexponential volume growth  of $M$ (due to Assumption $(C)$) to conclude that
$$
\int_M u(x, 1) d\mu(x) = \int_M p^V_1(x, y) d\mu(x) \le C.
$$The rest of the proof is identical.

\medskip

{\it Case 2.}  Assume $V \in L^p(M, \mu)$ with $1<p<2$. 

This is almost identical to that of Case 2, Proposition \ref{1,1}. Indeed, from that case, we have
$$\int_ M u(x, t)
d\mu(x) \le C + C' \sqrt{ t } \left(\int^t_1 \int_ M u^2(x, s) d\mu(x)ds  \right)^{(p-1)/2}
$$Since $p^V_t$ is contractive in $L^2$, this implies
$$
\int_ M u(x, t)
d\mu(x) \le C t^{p/2}, \quad t \ge 1.
$$

\medskip
{\it Case 3.}  Assume $V \in L^p(M, \mu)$ with $p \ge 2$.

Let $u=u(x, t)$ be as above. Then, as before, 
\[
\aligned
  \int_ M
u(x, t) d\mu(x) &\le  \int_ M u(x, 1) d\mu(x)
+ \int^t_1\int_ M V(x)
u(x, s) d\mu(x)ds\\
& \le C 
+ \int^t_1 \Vert V \Vert_p \
\Vert u(\cdot, s) \Vert_{p/(p-1)}ds\\
&\le C +   C\Vert V
\Vert_p \int^t_1
\Vert  P^V_s \Vert_{1,
p/(p-1)} ds.
\endaligned
\]Applying Riesz-Thorin interpolation  with the
parameters
\[
p_2=1, q_2=p/(p-1); \  p_0=1, q_0=2; \  p_1=1,
q_1=1; \
\theta=(2-q_2)/q_2,
\]we have
\[
\frac{1}{p_2}=
\frac{1-\theta}{p_0} +  \frac{\theta}{p_1},
\frac{1}{q_2}=
\frac{1-\theta}{q_0} +  \frac{\theta}{q_1}
\]and
\[
\Vert  P^V_s \Vert_{1, p/(p-1)}  = \Vert  P^V_s \Vert_{p_2,
q_2} \le \Vert  P^V_s \Vert^{1-\theta}_{p_0, q_0} \ \Vert
 P^V_s \Vert^{\theta}_{p_1, q_1} = \Vert  P^V_s
\Vert^{2/p}_{1, 2} \Vert  P^V_s \Vert^{1-(2/p)}_{1, 1}.
\]Therefore
\[
\int_ M |u(x, t)| d\mu(x)  \le C  +\Vert V
\Vert_p
\int^t_1
\Vert  P^V_s \Vert^{2/p}_{1, 2} \ \Vert
 P^V_s
\Vert^{1-(2/p)}_{1, 1}  ds 
\]Notice that
\[
\Vert  P^V_s \Vert_{1, 2} \le \Vert  P^V_s
\Vert_{2, \infty}
\le \Vert  P^V_{s-1} \Vert_{2, 2} \ \Vert  P^V_1
\Vert_{2,
\infty} \le C.
\]We obtain
\[
\int_ M |u(x, t)| d\mu(x)  \le C
 + \Vert V
\Vert_p\int^t_1
   \Vert  P^V_s
\Vert^{1-(2/p)}_{1, 1}  ds.
\]i.e.
\[
\Vert  P^V_t \Vert_{1, 1}
\le C + \Vert V
\Vert_p\int^t_1 \Vert  P^V_s
\Vert^{1-(2/p)}_{1, 1}  ds.
\]From here it is easy to see that
\[
\Vert  P^V_t \Vert_{1, 1}
\le C t^{p/2}.
\]This completes Step 1.

\bigskip
{\it Step 2. }
Write
\[
I(t) = \int_M u^2(x, t) \,d\mu(x).
\]
As in the proof of Proposition \ref{nash}, the strong positivity of $\Delta-V$ yields
\begin{equation}
I(t) \le
\frac{\int_{ \{x | u(x, t) >s \}} |\nabla (u(x,
t)-s)|^2
\,d\mu(x)}{\lambda_1 ( \{x | u(x, t) >s \})}  + 2s
F(t).
\label{lev1}
\end{equation}Here, according to \eqref{1to1}, $F(t) =t^{p/2}$.

Using
the fact that
\[
| \{x | u(x, t) >s \} | \le s^{-1} \int_M u (x, t)
d\mu(x),
\]and $(UFK)$, we deduce
\[
I(t)
\le \frac{\int_{ \{x | u(x, t) >s \}} |\nabla (u(x,
t)-s)|^2
\,d\mu(x)}{\Lambda (s^{-1} F(t))}  + 2s
F(t).
\]

Hence
\begin{equation}
\int_{ \{x | u(x, t) >s \}} |\nabla u(x,
t)|^2 \,d\mu(x) \ge [
I(t) - 2 s F(t) ] \Lambda (s^{-1} F(t)) .
\label{sub1}
\end{equation}

By the strong positivity of $\Delta-V$, we have as usual
$$I'(t) \le -2 (1-A) \int_ M |\nabla u(x, t)|^2 \,d\mu(x),$$
thus  the combination of  the above inequalities
yields
\begin{equation}
I'(t) \le - 2(1- A) [ I(t) - 2 s F(t) ] \Lambda
(s^{-1} F(t)).
\label{comb}
\end{equation}

Take $s F(t) = I(t)/4$,
i.e.
\[
s^{-1} = 4 I^{-1}(t) F(t).
\] Then (\ref{comb})
becomes
\begin{equation}
I'(t) \le - (1- A) I(t) \Lambda (4 F^2(t)
I^{-1}(t)). \label{bec1}
\end{equation}
Hence
\[
\int^{2t}_t \frac{
I'(l)}{I(l) \Lambda(4 F^2(l) I^{-1}(l))} dl \le
-(1- A) t.
\]Notice that
$\Lambda$ is a decreasing and $F$ is an increasing
function. Therefore, for
$l \ge t$,
\[
\Lambda (4 F^2(l) I^{-1}(l)) \le \Lambda (4 F^2(t)
I^{-1}(l))
\]Consequently
\[
\int^{2t}_t \frac{ I'(l)}{I(l) \Lambda(4 F^2(t)
I^{-1}(l))} dl \le
-(1- A) t.
\]Take $\eta = 4 F^2(t) I^{-1}(l)$. One
gets
\[
\int^{ 4 F^2(t) I^{-1}(2t)}_{4 F^2(t) I^{-1}(t)} \frac{
d\eta}{\eta
\Lambda(\eta)}  \ge (1- A) t.
\]Following the definition of $\gamma$,
i.e. $t = \int^{\gamma}_0 \frac{
d\eta}{\eta \Lambda(\eta)} $, we
have
\[
\frac{ 4 F^2(t)}{I(2 t)} \ge \gamma ((1- A) t),
\]i.e.
\[
I(t) \le
\frac{4 F^2(t)}{\gamma( c t)}.
\]From here the 
desired bound for $\vec{p_t}$ follows immediately. \qed

\bigskip

Let us conclude by writing a semigroup version of the last part of the proof of Theorem \ref{final}, in the spirit of \cite{CUl}, where the case $F$ bounded is treated. We leave the details to the reader.

\begin{proposition} Let $(M,\mu)$ a $\sigma$-finite measure space, and $T_t$ be a semigroup acting  on $L^p(M,\mu)$,  for $1\leq p\leq
+\infty$,
  with infinitesimal generator $-A$.  Suppose that  there exists a non-decreasing function $F$ on $\reel_+$ such that
  $$\|T_t\|_{1\to 1}, \|T_t\|_{\infty\to \infty}\leq
F(t),\ \forall\,t>0,$$ and that
$$\theta(\|f\|_2^2)\leq {\rm Re}(Af,f),\,  \forall f\in D(A),\,
\|f\|_1\leq C,$$
fro some $C>0$,
where $\theta:]0,+\infty[\to]0,+\infty[$ is continuous and  satisfies
$\int^{+\infty}\frac{dx}{ \theta(x)}<+\infty$.
Then $T_t$ is ultracontractive and
$$\|T_t\|_{1\to \infty}\leq CF^2(t)m (Ct),\,\forall t>0,$$
for some $C>0$, 
where $m$ is the solution of $$-m'(t)=\theta(m(t))$$ on $]0,+\infty[$
such that $m(0)=+\infty$, or alternatively the inverse function
of $p(t)=\int_t^{+\infty}\frac{dx}{ \theta(x)}$.
\end{proposition}

\bigskip

{\bf Acknowledgement:} The second author acknowledges the support of the University of Cergy-Pontoise during the preparation of this paper. Both authors thanks Adam Sikora for nice remarks on the manuscript.

\enddocument